\documentclass[a4paper,11pt,reqno]{amsart}
\usepackage[T1]{fontenc}
\usepackage{lmodern}
\usepackage{amsmath,amssymb,amsfonts,amsthm,bm,mathtools}
\usepackage{enumitem}
\usepackage[final]{microtype}
\microtypesetup{nopatch=footnote}
\usepackage{xurl}
\usepackage[numbers,sort&compress]{natbib}
\usepackage[hidelinks]{hyperref}
\hypersetup{pdftitle={Relative computability and definability for nonregular oracles},pdfauthor={Yong Liu, Cheng Peng, and Mengzhou Sun}}
\numberwithin{equation}{section}
\newtheorem{theorem}{Theorem}[section]
\newtheorem{lemma}[theorem]{Lemma}

\newtheorem{claim}{Claim}

\theoremstyle{definition}
\newtheorem{definition}[theorem]{Definition}
\theoremstyle{remark}

\newcommand{\tup}[1]{\langle #1 \rangle}
\newcommand{\sT}{\operatorname{sT}}
\newcommand{\pT}{\operatorname{pT}}
\newcommand{\CEpm}{\operatorname{CE}^{\pm}}
\DeclareMathOperator{\dom}{dom}
\newcommand{\res}{\!\!\restriction\!}
\DeclareMathOperator{\cf}{cf}
\newcommand{\abs}[1]{\left| #1 \right|}
\newcommand{\ext}{\mathop{\mathrm{end}}}
\DeclareMathOperator{\bit}{bit}
\newcommand{\cU}{\mathcal{U}}
\newcommand{\cC}{\mathcal{C}}

\begin{document}

\title{Relative computability and definability for nonregular oracles}
\author[Liu]{Yong Liu}
\address[Liu]{School of Information Engineering\\
Nanjing Xiaozhuang University\\
CHINA}
\email{liuyong@njxzc.edu.cn}
\author[Peng]{Cheng Peng}
\address[Peng]{Institute of Mathematics\\
Hebei University of Technology\\
CHINA}
\email{pengcheng@hebut.edu.cn}
\author[Sun]{Mengzhou Sun}
\address[Sun]{Institute of Mathematics, University of Warsaw, Poland}
\email{m.sun3@uw.edu.pl}
% \thanks{Liu's research was partially funded by Nanjing Xiaozhuang University (No.~2022NXY39). Peng's research was partially funded by the National Natural Science Foundation of China (No.~12671001). Sun's work was supported by grant no. 2023/49/B/ST1/02627 of the National Science Centre, Poland.}

\begin{abstract}
Classically, a set is computable from an oracle if both the set and its
complement are computably enumerable (c.e.) relative to that oracle. We
construct a countable nonstandard model of true arithmetic with sets \(A\)
and \(X\) such that both \(A\) and its complement are c.e.\ relative to
\(X\), yet \(A\) is not pointwise Turing reducible to \(X\). Here
we use Turing functionals that are single-valued at every oracle. For
regular oracles, we prove that the classical equivalence holds over
\(I\Sigma_1\). Two further phenomena, each known from earlier work, can be
realized with the same oracle \(X\) that witnesses our main separation.
Strong and pointwise Turing reducibility differ for this \(X\), and some
\(\Delta_1^X\) set is not \(X\)-c.e.
\end{abstract}
\maketitle
\enlargethispage{2pt}

\section{Introduction}\label{sec:intro}

In classical computability theory, \(A\) is computable relative to
\(X\) when a Turing functional computes the characteristic function
of \(A\) from \(X\). It is known that this holds if and only if
both \(A\) and its complement are computably enumerable (c.e.)
relative to \(X\): given an input, we run the two enumerations until
it appears in one of them. However, over a nonstandard model of arithmetic, the validity
of this argument depends on the definition of relative computability
and on the properties of the oracle.

There are several related notions in the literature. Mytilinaios and
Slaman~\cite[Definition~3.2]{RN291} define relative recursiveness for
sets by the existence of the two relative enumerations. Pointwise
Turing reducibility is instead defined using a single c.e.\ set of
axioms. An axiom specifies an input, an output, and two \(M\)-finite sets
of oracle queries, one requiring positive answers and the other negative
answers. We adopt the monotonicity
and single-valuedness conditions of Groszek and Hummel~\cite[\S2.1]{RN323},
also used in~\cite{RN419}. These conditions make a Turing functional
single-valued at \emph{every} oracle. If single-valuedness is required only at the
particular oracle under consideration, as in the convention
of~\cite[Definition~2.8]{RN324}, the two relative enumerations can
always be combined. Under the conditions we adopt, this need not be
possible.

For an oracle \(X\), write \(\pT(X)\) for the sets pointwise Turing
reducible to \(X\), and \(\CEpm(X)\) for the sets whose members and
nonmembers can both be enumerated relative to \(X\). If \(X\) is
regular, meaning that the intersection of \(X\) with each bounded initial
segment of the model is coded, then \(\pT(X)=\CEpm(X)\) by Theorem~\ref{thm:regular}.
We write \(\sT(X)\) for the sets strongly Turing reducible to \(X\),
and \(\Delta_1^X\) for the sets \(A\) such that both \(A\) and its
complement are \(\Sigma_1^X\). The precise definitions are in Section~\ref{sec:prelim}.

\begin{theorem}\label{thm:main}
There exists a countable model \(N\equiv\mathbb N\) and \(X\subseteq N\)
such that all three of the following hold.
\begin{enumerate}[label=\textup{(\alph*)},ref=\alph*]
\item\label{it:main-strong} \(\sT(X)\subsetneq\pT(X)\).
\item\label{it:main-pointwise} \(\pT(X)\subsetneq\CEpm(X)\).
\item\label{it:main-definability} \(\CEpm(X)\subsetneq\Delta_1^X\).
\end{enumerate}
\end{theorem}

Part~(\ref{it:main-pointwise}) is the main point of Theorem~\ref{thm:main}.
For the oracle \(X\) we construct, a set \(A\) and its complement are
both \(X\)-c.e., yet no Turing functional computes the characteristic
function of \(A\) from \(X\). The other two separations were known separately. For part~(\ref{it:main-strong}), strong Turing
reducibility is transitive~\cite[Proposition~2.9]{RN324}, whereas pointwise Turing
reducibility can fail to be transitive~\cite{RN324}. For part~(\ref{it:main-definability}),
Michalski~\cite[Theorem~1.2]{RN424} gives a model of Peano arithmetic
with an oracle \(Y\) and a \(\Delta_1^Y\) set that is not \(Y\)-c.e.
Here both separations hold in the same countable model and for the
same oracle \(X\) as in part~(\ref{it:main-pointwise}).

Section~\ref{sec:prelim} gives the definitions and proves the equalities
for regular oracles. In Section~\ref{sec:main-proof} we first prove
Theorem~\ref{thm:main}(\ref{it:main-pointwise}) by constructing a column
oracle and passing to a countable model. We then show that the oracle can be chosen to satisfy
parts~(\ref{it:main-strong}) and~(\ref{it:main-definability}) as well.

\section{Preliminaries}\label{sec:prelim}

\subsection{Coding and relative formulas}
Unless otherwise stated, \(M\models I\Sigma_1\). We use \(\mathbb N\)
to denote the standard model of arithmetic and \(\omega\) to denote
the external set of standard natural numbers. All cardinalities in
this paper are external. A set is \emph{externally finite} if its
cardinality is less than \(\aleph_0\), and \emph{externally infinite}
otherwise. We allow number parameters in definitions. For
\(A\subseteq M\), write \(\overline{A}=M\setminus A\), and identify \(A\)
with its characteristic function when it appears as the value of a
computation.

A set \(A\subseteq M\) is \emph{coded} if there exists a code \(c\in M\)
such that
\[
 x\in A\iff M\models\bit(c,x)=1,
\]
where \(\bit(c,x)\) is the \(x\)-th bit of \(c\). This relation is
expressible by a bounded formula, that is, a \(\Delta_0\) formula. A set is
\emph{\(M\)-finite} if it is coded. Such a set is bounded in \(M\),
but need not be externally finite. We identify coded sets with their
codes when quantifying over them. Finite unions and bounded arithmetic subsets
of coded sets are again coded.

Write \([0,a)^M=\{x\in M\mid M\models x<a\}\). An oracle is an
arbitrary subset \(X\subseteq M\). It is \emph{regular}, also called
\emph{amenable} or \emph{piecewise coded}, if \(X\cap[0,a)^M\) is \(M\)-finite for every
\(a\in M\).

A set \(A\subseteq M\) is \(\Sigma_1^X\) if, for some bounded formula
\(\theta\) in the language with a unary predicate for \(X\),
\[
 x\in A\iff(M,X)\models\exists u\,\theta(x,u).
\]
A c.e.\ set in \(M\) is defined to be an arithmetically
\(\Sigma_1\)-definable set with parameters. The same pattern is not
used for oracles: being \(\Sigma_1^X\) is not the definition of being
\(X\)-c.e. That definition is given below, in terms of Turing
functionals. The two notions agree when \(X\) is regular
(Theorem~\ref{thm:regular}(\ref{it:regular-formulas})). In particular
they agree for \(X=\varnothing\), since the empty set is regular.
The assumption \(M\models I\Sigma_1\) concerns
the arithmetic language; induction or collection with the oracle
predicate is not assumed.

\subsection{Turing functionals and relative computability}

Following~\cite{RN323,RN419}, we introduce the notion of Turing
functional.
\begin{definition}[Turing functional]\label{def:turing_functional}
A \emph{Turing functional} is a c.e.\ set \(\Phi\) of \emph{axioms}
\(\tup{x,y,P,N}\), where \(x,y\in M\), \(P,N\) are \(M\)-finite
sets, and the following conditions hold:
\begin{enumerate}[label=\textup{(F\arabic*)},ref=(F\arabic*),leftmargin=3.2em]
\item\label{it:disjoint}
\emph{Disjointness.}
\(\tup{x,y,P,N}\in\Phi\) implies \(P\cap N=\varnothing\).
\item\label{it:extension}
\emph{Monotonicity.}
If \(\tup{x,y,P,N}\in\Phi\) and \(P',N'\) are \(M\)-finite sets
with \(P\subseteq P'\), \(N\subseteq N'\), and
\(P'\cap N'=\varnothing\), then \(\tup{x,y,P',N'}\in\Phi\).
\item\label{it:uniqueness}
\emph{Single-valuedness.}
If \(\tup{x,y,P,N},\tup{x,y',P,N}\in\Phi\), then \(y=y'\).
\end{enumerate}
\end{definition}
Note that single-valuedness~\ref{it:uniqueness} is not required
in~\cite[Definition~2.8]{RN324}.

\begin{definition}\label{def:compatible_axioms}
Two axioms \(\tup{x,y,P,N}\) and \(\tup{x,y',P',N'}\) are \emph{compatible} if \((P\cup P')\cap(N\cup N')=\varnothing\), or equivalently, \(P\cap N'=\varnothing\) and \(P'\cap N=\varnothing\); otherwise, we call them \emph{incompatible}.
\end{definition}

\begin{lemma}\label{lem:opposite-outputs}
Let \(\Phi\) be a Turing functional. If \(\tup{x,y,P,N},\tup{x,y',P',N'}\in\Phi\) are compatible, then \(y=y'\).
\end{lemma}
\begin{proof}
Since \((P\cup P')\cap(N\cup N')=\varnothing\), we have \(P\cap N'=\varnothing\) and \(P'\cap N=\varnothing\). By monotonicity~\ref{it:extension}, we have \(\tup{x,y,P\cup P',N\cup N'}\in\Phi\) and \(\tup{x,y',P\cup P',N\cup N'}\in\Phi\). Single-valuedness~\ref{it:uniqueness} implies \(y=y'\).
\end{proof}

We write \(\Phi^X(x)=y\) if 
\begin{equation}\label{eq:computation}
 \exists P,N\,
 \bigl(\tup{x,y,P,N}\in\Phi\ \land\ P\subseteq X
                         \ \land\ N\cap X=\varnothing\bigr).
\end{equation}
We say \(\Phi^X(x)\downarrow\) if \(\exists y\,(\Phi^X(x)=y)\); \(\Phi^X(x)\uparrow\) if \(\neg(\Phi^X(x)\downarrow)\).
We write \(\dom\Phi^X=\{x\mid\Phi^X(x)\downarrow\}\).
The notation \(\Phi^X(x)\ne y\) means \(\neg(\Phi^X(x)=y)\), equivalently, \(\Phi^X(x)\uparrow\) or \(\Phi^X(x)\downarrow\neq y\).
A set \(A\) is \(X\)-c.e.\ if \(A=\dom\Phi^X\) for some Turing
functional \(\Phi\). Clearly, every \(X\)-c.e.\ set is \(\Sigma_1^X\).
However, a \(\Sigma_1^X\) set need not be \(X\)-c.e.~\cite{RN424}.

Fix a uniform c.e.\ enumeration \(\{\Phi_e\mid e\in M\}\) of all Turing
functionals over \(M\), possibly with repetitions.

\begin{definition}\label{def:turing_reducibility}
Let \(A,B\subseteq M\).
\begin{enumerate}[label=\textup{(\alph*)},ref=\alph*]
\item \(A\) is \emph{pointwise Turing reducible} to \(B\), denoted by
\(A\leq_{pT}B\), if there is a Turing functional \(\Phi\) such that
\(A=\Phi^B\). That is, \(A(x)=\Phi^B(x)\) for all \(x\in M\).
\item \(A\) is \emph{strongly Turing reducible} to \(B\), denoted by
\(A\leq_{sT}B\), if the sets of codes
\[
 \begin{aligned}
 A^+&=\{P\mid P\text{ is an }M\text{-finite subset of }A\},\\
 A^-&=\{N\mid N\text{ is an }M\text{-finite subset of }\overline{A}\}
 \end{aligned}
\]
are pointwise Turing reducible to \(B\).
\end{enumerate}
\end{definition}

For \(X\subseteq M\), define
\begin{align*}
 \sT(X)&=\{A\subseteq M\mid A\leq_{sT}X\},\\
 \pT(X)&=\{A\subseteq M\mid A\leq_{pT}X\},\\
 \CEpm(X)&=\{A\subseteq M\mid A\text{ and }\overline{A}\text{ are both }X\text{-c.e.}\},\\
 \Delta_1^X&=\{A\subseteq M\mid A\text{ and }\overline{A}\text{ are both }\Sigma_1^X\}.
\end{align*}
The model under consideration is understood from context and is omitted
from these notations. It follows from the definitions that
\begin{equation}\label{eq:inclusions}
 \sT(X)\subseteq\pT(X)\subseteq\CEpm(X)\subseteq\Delta_1^X.
\end{equation}
We do not define what \(X\)-computable means in this paper, as any of
these classes is a candidate.\footnote{In~\cite{RN291,RN303}, a set \(A\) is \(X\)-computable if \(A\in\CEpm(X)\). In~\cite{RN419}, \(A\) is \(X\)-computable if \(A\in \Delta_1^X\),
and the authors identify this condition with \(A\in \pT(X)\).
For regular oracles, including the c.e.\ oracles over \(I\Sigma_1\)
used there, Theorem~\ref{thm:regular} gives
\(\pT(X)=\CEpm(X)=\Delta_1^X\). Theorem~\ref{thm:main} shows that
all three classes can be distinct for the same oracle.}
If \(X\) is regular and
\((M,X)\models B\Sigma_1^X\),
then \(\sT(X)=\pT(X)\); see~\cite[Proposition~2.9]{RN324}.
Michalski~\cite[p.~517, following (1.1)]{RN424} observes that regularity
suffices for part~(\ref{it:regular-formulas}) of the next theorem. For
completeness, we give proofs of both parts.

\subsection{The regular case}

\begin{theorem}\label{thm:regular}
Let \(X\subseteq M\) and \(M\models I\Sigma_1\).
\begin{enumerate}[label=\textup{(\alph*)},ref=\alph*]
\item\label{it:regular-merge}
If \(X\) is regular, then \(\pT(X)=\CEpm(X)\).
\item\label{it:regular-formulas}
If \(X\) is regular, then every \(\Sigma_1^X\) set is \(X\)-c.e.,
and therefore \(\CEpm(X)=\Delta_1^X\).
\end{enumerate}
\end{theorem}

\begin{proof}
\textbf{(a)} Let \(A=\dom\Phi^X\) and \(\overline{A}=\dom\Psi^X\).
Given a binary string \(\sigma\) of length \(s\), we define \[
    P_\sigma=\{z<s\mid \sigma(z)=1\} \text{ and } N_\sigma=\{z<s\mid \sigma(z)=0\}.
\]
We write \(\sigma\preccurlyeq\tau\) to denote that \(\sigma\) is an initial segment of \(\tau\). An axiom of the form \(\tup{x,y,P,N}\) is called an \(x\)-axiom.

We first construct a c.e.\ set \(\Gamma^*\) of axioms of the form
\[
    \tup{x,y,P_\sigma,N_\sigma}
\]
such that any two distinct \(x\)-axioms in \(\Gamma^*\) have incompatible
strings. That is, if their strings are \(\sigma\) and \(\tau\), then
\(\sigma(k)\neq\tau(k)\) for some \(k<\min\{|\sigma|,|\tau|\}\).
Equivalently, these axioms are incompatible in the sense of
Definition~\ref{def:compatible_axioms}. The \emph{closure} of \(\Gamma^*\)
consists of all \(\tup{x,y,P,N}\) such that, for some \(\sigma\),
\[
    \tup{x,y,P_\sigma,N_\sigma}\in \Gamma^* \text{ and } P_\sigma\subseteq P \text{ and }N_\sigma\subseteq N\text{ and } P\cap N=\varnothing.
\]

Let \(\{\Phi_s\}_{s\in M}\) and \(\{\Psi_s\}_{s\in M}\) be the enumerations of
\(\Phi\) and \(\Psi\) respectively. 

At stage \(s\), for each \(x<s\), consider each coded string \(\sigma\)
of length \(s\) that does not extend any previously \(x\)-used string.
An \(x\)-axiom \(\tup{x,y,P,N}\) applies to \(\sigma\) if
\(P\subseteq P_\sigma\) and \(N\subseteq N_\sigma\).
\begin{enumerate}
    \item\label{it:add-axiom1}
    If some \(x\)-axiom in \(\Phi_s\) applies to \(\sigma\), enumerate
    \(\tup{x,1,P_\sigma,N_\sigma}\) into \(\Gamma^*\).
    \item\label{it:add-axiom0}
    Otherwise, if some \(x\)-axiom in \(\Psi_s\) applies to \(\sigma\),
    enumerate \(\tup{x,0,P_\sigma,N_\sigma}\) into \(\Gamma^*\).
\end{enumerate}
If either rule enumerates an axiom, we say that \(\sigma\) is
\(x\)-\emph{used}. Distinct \(x\)-used strings are incompatible:
strings used at the same stage have the same length, and a string
used at a later stage cannot extend one used earlier. Thus the
closure \(\Gamma\) of \(\Gamma^*\) is a Turing functional.

Next, we show that \(A=\Gamma^X\). Fix \(x\in M\) and suppose
\(A(x)=1\). Choose \(s>x\) large enough that some axiom
\(\tup{x,y,P,N}\in\Phi_s\) witnesses \(\Phi^X(x)\downarrow\) and
\(P\cup N\subseteq[0,s)^M\). By regularity, \(\sigma=X\res s\) is
coded, and this axiom applies to \(\sigma\). 
Note that for any \(\tau\preccurlyeq \sigma\), no \(x\)-axiom in \(\Psi\) applies to \(\tau\), since \(\Psi^X(x)\uparrow\). Therefore, there must exist some \(\tau\preccurlyeq\sigma\) such that Item~(\ref{it:add-axiom1}) happens and therefore \(\tup{x,1,P_\tau,N_\tau}\in \Gamma\) and \(\Gamma^X(x)=1\).

\smallskip
\noindent\textbf{(b)} Say \(x\in A\) if and only if
\(\exists u\,\theta(x,u)\), where \(\theta\) is a bounded formula in
the language with the oracle predicate.
Put \(\theta\) in negation normal form, so that negations occur only
at atomic formulas. For disjoint \(M\)-finite sets \(P,N\), let
\(\theta^{P,N}(x,u)\) be obtained by replacing each oracle literal
\(t\in X\) by \(t\in P\), and each \(t\notin X\) by \(t\in N\).
This is a bounded arithmetic formula and has the following properties:
\begin{enumerate}[label=\textup{\alph*.},leftmargin=2em]
\item \emph{Monotonicity.} If \(P\subseteq P'\), \(N\subseteq N'\),
and \(P'\cap N'=\varnothing\), then
\(\theta^{P,N}(x,u)\) implies \(\theta^{P',N'}(x,u)\).
\item \emph{Soundness.} If \(P\subseteq X\) and
\(N\cap X=\varnothing\), then \(\theta^{P,N}(x,u)\) implies
\(\theta(x,u)\).
\end{enumerate}
Enumerate a Turing functional \(\Theta\) by putting
\begin{equation}\label{eq:formula-axioms}
 \tup{x,1,P,N}\in\Theta
 \iff P\cap N=\varnothing\ \land\ \exists u\,\theta^{P,N}(x,u).
\end{equation}
By the monotonicity of \(\theta^{P,N}\), the set \(\Theta\) satisfies
monotonicity~\ref{it:extension}, and single-valuedness holds since the
only output is \(1\).

We show that \(A=\dom\Theta^X\).
For \(x\in A\), choose \(u\) with \(\theta(x,u)\). The formula
itself gives an arithmetic bound \(d\) such that its truth depends
only on the first \(d\) bits of \(X\): bound the arguments of all
oracle atoms over the bounded quantifiers of \(\theta\). Let
\(P=X\cap[0,d)^M\) and \(N=[0,d)^M\setminus P\). Regularity makes
both sets \(M\)-finite, and \(\theta^{P,N}(x,u)\) holds. By
\eqref{eq:formula-axioms}, \(x\in\dom\Theta^X\).
Conversely, an axiom witnessing \(x\in\dom\Theta^X\) gives
\(\theta^{P,N}(x,u)\) for some \(u\), with \(P\subseteq X\) and
\(N\cap X=\varnothing\). Soundness gives \(\theta(x,u)\), so
\(x\in A\). Applying this also to \(\overline{A}\) proves
\(\CEpm(X)=\Delta_1^X\).
\end{proof}

If we omit single-valuedness~\ref{it:uniqueness} in
Definition~\ref{def:turing_functional}, as in~\cite{RN324}, then
\(\pT(X)=\CEpm(X)\) for every oracle. Indeed, if
\(A=\dom\Phi^X\) and \(\overline{A}=\dom\Psi^X\), put
\begin{align*}
 \tup{x,1,P,N}\in\Gamma&\iff\exists y\,\tup{x,y,P,N}\in\Phi,\\
 \tup{x,0,P,N}\in\Gamma&\iff\exists y\,\tup{x,y,P,N}\in\Psi.
\end{align*}
Then \(\Gamma^X=A\), although \(\Gamma^Y\) may be multivalued for
another oracle \(Y\).

\section{Proof of the main theorem}\label{sec:main-proof}

We first construct a model \(M\equiv\mathbb N\) whose cofinality
allows us to take unions of bounded configurations at limit stages.
In this model we construct a column oracle \(X\) and a set
\(A_X\in\CEpm(X)\setminus\pT(X)\). Once we have this separation in
a countable elementary substructure, we return to the construction
in \(M\) to obtain all three separations.

\subsection{A \texorpdfstring{\(\kappa^+\)}{κ⁺}-like model for the oracle construction}\label{sec:model}

For a cardinal \(\kappa\), write \(\kappa^+\) for its successor cardinal.
The \emph{cofinality} of \(M\), denoted by \(\cf(M)\), is the least
cardinality of a cofinal subset of \(M\). Thus every subset of
cardinality less than \(\cf(M)\) is bounded in \(M\). We write
\(M_0\prec_{\ext}M\) if \(M\) is a proper \emph{elementary end-extension}
of \(M_0\): \(M\setminus M_0\neq \varnothing\) and \(M_0\prec M\) and \(M_0\) is an initial segment of
\(M\). The initial segment condition means that every element of
\(M\setminus M_0\) is greater than every element of \(M_0\).

\begin{definition}\label{def:like}
Let \(\lambda\) be an infinite cardinal. A model \(M\) is
\emph{\(\lambda\)-like} if \(\abs M=\lambda\) and
\(\abs{[0,a)^M}<\lambda\) for every \(a\in M\); see~\cite{Kaye1997}.
\end{definition}

The following lemma gives a \(\kappa^+\)-like model with the additional
properties needed in the construction.
\begin{lemma}\label{lem:main1}
Let \(\kappa=2^{\aleph_0}\). There exists a \(\kappa^+\)-like model
\(M\) such that
\begin{enumerate}[label=\textup{(L\arabic*)},ref=(L\arabic*),leftmargin=3.2em]
\item\label{it:like1} \(M\equiv\mathbb N\);
\item\label{it:like2} for every nonstandard \(a\in M\),
\(\abs{[0,a)^M}=\kappa\);
\item\label{it:like3} \(\abs M=\cf(M)=\kappa^+\);
\item\label{it:like4} every externally infinite \(M\)-finite subset of
\(M\) has cardinality \(\kappa\).
\end{enumerate}
\end{lemma}

\begin{proof}
We first note that~\ref{it:like3} follows from \(\kappa^+\)-likeness.
Indeed, suppose \(C\subseteq M\) were cofinal with \(\abs C<\kappa^+\). Then
\[
 M=\bigcup_{a\in C}[0,a)^M.
\]
Each initial segment on the right has cardinality at most \(\kappa\), and
there are at most \(\kappa\) of them, so the union has cardinality at
most \(\kappa\). This contradicts \(\abs M=\kappa^+\).
Condition~\ref{it:like4} follows from~\ref{it:like1} and
\ref{it:like2}. To see this, let \(E\) be an externally infinite
\(M\)-finite set. Inside \(M\), its elements have an increasing
enumeration with domain \([0,a)^M\), where \(a\) is its internal
cardinality. This is also a bijection externally. Since \(E\) is
externally infinite, \(a\) is nonstandard, and therefore
\(\abs E=\abs{[0,a)^M}=\kappa\).
Thus it suffices to construct a \(\kappa^+\)-like model satisfying
\ref{it:like1} and~\ref{it:like2}.

\smallskip
\noindent\textbf{Step 1. Start with a nonstandard model.}
Let \(\cU\) be a nonprincipal ultrafilter on \(\omega\) and take the
ultrapower of \(\mathbb N\)
\[
 M_0=\mathbb N^\omega/\cU.
\]
An element of \(\mathbb N^\omega\) is a function from \(\omega\) to
\(\mathbb N\), denoted by \(\lambda i.f(i)\). Two functions
\(\lambda i.f(i)\) and \(\lambda i.g(i)\) are equivalent modulo
\(\cU\) if
\[
 \{i\in\omega\mid f(i)=g(i)\}\in\cU.
\]
The elements of \(M_0\) are the equivalence classes
\([\lambda i.f(i)]_\cU\). The ordering is given by
\[
 [\lambda i.f(i)]_\cU<[\lambda i.g(i)]_\cU
 \iff\{i\in\omega\mid f(i)<g(i)\}\in\cU.
\]
We identify \(n\in\mathbb N\) with the class of the constant
function \([\lambda i.n]_\cU\). By \L o\'s's theorem this identifies
\(\mathbb N\) with an elementary substructure of \(M_0\), and in
fact \(\mathbb N\prec_{\ext}M_0\). Also,
\([\lambda i.n]_\cU<[\lambda i.i]_\cU\) for every standard \(n\),
so \(M_0\) is nonstandard. Clearly, \(\abs{M_0}\leq\kappa\).

\begin{claim}\label{claim:counting}
For every nonstandard \(a\in M_0\), \(\abs{[0,a)^{M_0}}=\kappa\).
\end{claim}
\begin{proof}[Proof of Claim~\ref{claim:counting}]
Let \(a=[\lambda i.a_i]_\cU\) be nonstandard. Define
\[
 c:(0,1)\longrightarrow[0,a)^{M_0},\qquad
 c(r)=[\lambda i.\lfloor r a_i\rfloor]_\cU,
\]
where \((0,1)\) is the interval of real numbers and \(\lfloor r a_i\rfloor\)
is the greatest integer at most \(r a_i\). Since
\(\{i\mid a_i>0\}\in\cU\), we have \(c(r)<a\).
For \(0<r<s<1\), choose a standard \(K\) with \(rK+1<sK\).
The set \(V_K=\{i\mid a_i>K\}\) belongs to \(\cU\). For every
\(i\in V_K\), \(ra_i+1<sa_i\), and hence
\(\lfloor ra_i\rfloor<\lfloor sa_i\rfloor\). Thus \(c(r)<c(s)\)
and \(c\) is injective. Consequently,
\[
 \kappa\leq\abs{[0,a)^{M_0}}\leq\abs{M_0}\leq\kappa.\qedhere
\]
\end{proof}
In particular, \(\abs{M_0}=\kappa\).

\smallskip
\noindent\textbf{Step 2. End-extension.}
The MacDowell--Specker theorem states that every model of \(PA\) has
a proper elementary end-extension~\cite[Theorem~1.1]{Schmerl1973}.
We construct a continuous chain \(\{M_\alpha\mid\alpha<\kappa^+\}\)
such that
\[
 M_\alpha\prec_{\ext}M_{\alpha+1},\qquad\abs{M_\alpha}=\kappa,
\]
starting with \(M_0\) from Step 1. At a nonzero limit \(\beta<\kappa^+\),
let \(M_\beta=\bigcup_{\alpha<\beta}M_\alpha\). This is an elementary
end-extension of the earlier models, and
\(\abs{M_\beta}=\kappa\), since \(\abs\beta\leq\kappa\).

At a successor stage, take a proper elementary end-extension \(N\)
of \(M_\alpha\) and choose \(d_\alpha\in N\setminus M_\alpha\).
By downward L\"owenheim--Skolem, choose \(M_{\alpha+1}\prec N\)
of cardinality \(\kappa\) with
\(M_\alpha\cup\{d_\alpha\}\subseteq M_{\alpha+1}\).
Then \(M_\alpha\prec_{\ext}M_{\alpha+1}\).
Let
\[
 M=\bigcup_{\alpha<\kappa^+}M_\alpha.
\]
The upper bound on its cardinality is \(\kappa^+\), and the pairwise
distinct elements \(d_\alpha\) give the lower bound. Thus
\(\abs M=\kappa^+\) and
\(\mathbb N\prec_{\ext}M_0\prec_{\ext}M\).

\smallskip
\noindent\textbf{Step 3. Cardinality of initial segments.}
Let \(a\in M\) be nonstandard, and choose \(\alpha\) with
\(a\in M_\alpha\). Since \(M\) is an end-extension of \(M_\alpha\),
\[
 [0,a)^M=[0,a)^{M_\alpha},
\]
so its cardinality is at most \(\kappa\). If \(a\in M_0\), the
lower bound follows from Claim~\ref{claim:counting}. If
\(a\notin M_0\), then \(M_0\subseteq[0,a)^M\), which gives the
same lower bound. Thus \(M\) is \(\kappa^+\)-like and satisfies
\ref{it:like1} and~\ref{it:like2}, completing the construction.
\end{proof}

\subsection{Separating \texorpdfstring{\(\pT(X)\) from \(\CEpm(X)\)}{pT(X) from CE±(X)}}\label{sec:construction}

\begin{lemma}\label{lem:construction}
Let \(\kappa=2^{\aleph_0}\) and \(M\) be a \(\kappa^+\)-like
model satisfying conditions~\ref{it:like1}--\ref{it:like4} of
Lemma~\ref{lem:main1}. There exist \(A,X\subseteq M\) such that
\(A\in\CEpm(X)\setminus\pT(X)\).
\end{lemma}

The construction is carried out in ZFC. The idea is that if
\(A=\Phi^X\) for some Turing functional \(\Phi\), then
single-valuedness~\ref{it:uniqueness} will yield a contradiction. Fix a
nonstandard \(b\in M\) and let \(I=[0,b)^M\). We construct a set
\(D\subseteq I\), use it to define the column oracles, and then build
the desired oracle by extending configurations.

\subsubsection*{Construction of \(D\)}

\begin{lemma}\label{lem:D}
There exists \(D\subseteq I\) such that if \(E\subseteq I\) is
\(M\)-finite and either \(E\subseteq D\) or \(E\cap D=\varnothing\),
then \(E\) is externally finite. Moreover, \(D\) is not regular.
\end{lemma}
\begin{proof}
Every \(M\)-finite subset of \(I\) has a binary code below \(2^b\).
By Lemma~\ref{lem:main1}~\ref{it:like2},
\(\abs{[0,2^b)^M}=\kappa\). List all externally infinite
\(M\)-finite subsets of \(I\) as \(\{E_\alpha\}_{\alpha<\kappa}\),
allowing repetitions if there are not enough distinct ones. This family is nonempty since it contains \(I\).
At stage \(\alpha\), pick distinct elements
\[
 r_\alpha,s_\alpha\in
 E_\alpha\setminus\{r_\beta,s_\beta\mid\beta<\alpha\}.
\]
This is possible since \(\abs{E_\alpha}=\kappa\) by
Lemma~\ref{lem:main1}~\ref{it:like4}, whereas fewer than
\(\kappa\) elements have been chosen. Let
\(D=\{r_\alpha\mid\alpha<\kappa\}\). Each \(E_\alpha\) meets both
\(D\) and its complement, at \(r_\alpha\) and \(s_\alpha\),
respectively.

If \(D\) were regular, then \(D=D\cap[0,b)\) is a coded externally infinite subset of \(D\), a contradiction.
\end{proof}

\subsubsection*{Column oracles and the sets \(A_Y\)}

Use the Cantor pairing function
\(\tup{k,x}=(k+x)(k+x+1)/2+x\). It is a bijection in \(M\) and
is at least as large as either coordinate. Its graph has a bounded
arithmetic definition. For \(Y\subseteq M\), its \(k\)-th column is
\[
 Y^{[k]}=\{x\in M\mid\tup{k,x}\in Y\}.
\]
We say that \(Y\) is a \emph{column oracle} if
\begin{enumerate}[label=\textup{(\alph*)}]
\item \(Y^{[0]}=D\);
\item for each \(n\in M\), there exists \(k_n\in I\) such that
one of the following holds:
\begin{enumerate}[label=\textup{(\roman*)}]
\item \((Y^{[2n+1]},Y^{[2n+2]})=(D\Delta\{k_n\},D)\);
\item \((Y^{[2n+1]},Y^{[2n+2]})=(D,D\Delta\{k_n\})\).
\end{enumerate}
\end{enumerate}
Here \(D\Delta\{k_n\}\) denotes the symmetric difference. In
particular, every column is empty outside \(I\).

Note that a column oracle is not regular as its first column is \(D\), which is not \(M\)-finite.

For a column oracle \(Y\), define \(A_Y\) by
\begin{equation}\label{eq:A-Y}
\begin{aligned}
 n\in A_Y
 &\iff(Y^{[2n+1]},Y^{[2n+2]})=(D\Delta\{k_n\},D)\\
 &\iff\exists k<b\,\bigl(Y^{[0]}(k)\ne Y^{[2n+1]}(k)\bigr),\\
 n\notin A_Y
 &\iff(Y^{[2n+1]},Y^{[2n+2]})=(D,D\Delta\{k_n\})\\
 &\iff\exists k<b\,\bigl(Y^{[0]}(k)\ne Y^{[2n+2]}(k)\bigr).
\end{aligned}
\end{equation}
As usual, \(Y^{[j]}(k)\) denotes the value at \(k\) of the
characteristic function of the \(j\)-th column. For an \(M\)-finite set \(C\), we likewise write
\(C^{[j]}=\{x\in M\mid\tup{j,x}\in C\}\).

\begin{lemma}\label{lem:enumerations}
If \(Y\) is a column oracle, then \(A_Y\in\CEpm(Y)\).
\end{lemma}
\begin{proof}
Enumerate \(\tup{n,1,P,N}\) into \(\Gamma\) if
\(P\cap N=\varnothing\) and
\[
 \exists k<b\,
 \bigl(k\in P^{[0]}\cap N^{[2n+1]}
       \ \lor\ k\in N^{[0]}\cap P^{[2n+1]}\bigr).
\]
This is a c.e.\ set of axioms closed under consistent extensions.
Single-valuedness holds because it only outputs \(1\), so \(\Gamma\) is a
Turing functional. A disagreement between the two columns gives an
applicable axiom by taking \(P,N\) to be the two singleton query sets.
Thus \(A_Y=\dom\Gamma^Y\). Replacing \(2n+1\) by \(2n+2\) gives
a functional enumerating \(\overline{A_Y}\).
\end{proof}

\begin{lemma}\label{lem:externally-finite}
Suppose \(Y\) is a column oracle and \(\Gamma\) a Turing functional.
If \(\tup{n,i,P,N}\in\Gamma\) witnesses \(\Gamma^Y(n)=i\) and
\(Y^{[k]}=D\), then \(P^{[k]}\subseteq I\), and both \(P^{[k]}\) and
\(N^{[k]}\cap I\) are externally finite.
\end{lemma}
\begin{proof}
Since \(P\subseteq Y\), we have \(P^{[k]}\subseteq Y^{[k]}=D\subseteq I\),
and since \(N\cap Y=\varnothing\), we have
\((N^{[k]}\cap I)\cap D=\varnothing\). Both sets are \(M\)-finite by
bounded arithmetic comprehension, so Lemma~\ref{lem:D} applies.
\end{proof}

\subsubsection*{Configurations for column oracles}

We construct \(X\) by extending partial assignments. At each step we
fix some columns, and record this information in a configuration.
A \emph{configuration} is a partial function
\[
 \rho:M\rightharpoonup\{0,1\}\times I.
\]
We call \(\rho\) \emph{bounded} if its domain is bounded in \(M\).
Neither \(\rho\) nor its domain is required to have a code in \(M\).

For a configuration \(\rho\), let \(\cC(\rho)\) be the set of all
column oracles \(Y\) such that, for each \(n\in\dom\rho\),
\begin{enumerate}[label=\textup{(\alph*)}]
\item if \(\rho(n)=(1,k)\), then
\((Y^{[2n+1]},Y^{[2n+2]})=(D\Delta\{k\},D)\), so \(A_Y(n)=1\);
\item if \(\rho(n)=(0,k)\), then
\((Y^{[2n+1]},Y^{[2n+2]})=(D,D\Delta\{k\})\), so \(A_Y(n)=0\).
\end{enumerate}
The set \(\cC(\rho)\) is nonempty, since on each unassigned input we
can choose the column pair \((D\Delta\{0\},D)\). If \(\rho'\supseteq\rho\), then
\(\cC(\rho')\subseteq\cC(\rho)\).

The following is a ZFC lemma, not a lemma internal to \(M\).
\begin{lemma}[First extension lemma]\label{lem:extension}
Let \(\rho\) be a bounded configuration, \(\Gamma\) a Turing functional,
and \(n\in M\) with \(n\notin \dom \rho\). There exists a bounded configuration
\(\rho'\supseteq\rho\) with \(n\in \dom \rho'\) such that
\begin{equation}\label{eq:extension-conclusion}
 \Gamma^Z(n)\ne A_Z(n)\qquad\text{for every }Z\in\cC(\rho').
\end{equation}
Moreover, we can make \(\dom\rho'\) an initial segment
of \(M\).
\end{lemma}

\begin{proof}
Let \(\rho_0=\rho\). We try to
build a chain of bounded configurations
\(\rho_0\subseteq\rho_1\subseteq\cdots\), keeping \(n\notin\dom\rho_m\).

\smallskip
\noindent\textbf{Step 1. The question at stage \(m\).}
At a standard stage \(m\in\omega\), let
\(\rho_m^*=\rho_m\cup\{\tup{n,(1,m)}\}\). For every
\(Y\in\cC(\rho_m^*)\),
\[
 (Y^{[2n+1]},Y^{[2n+2]})=(D\Delta\{m\},D),
\]
and therefore \(A_Y(n)=1\). Ask whether some \(Y\in\cC(\rho_m^*)\) satisfies
\(\Gamma^Y(n)=1\).

\smallskip
\noindent\textbf{Step 2. The answer is negative.}
Suppose no such \(Y\) exists. Then \(\Gamma^Y(n)\ne1=A_Y(n)\) for
every \(Y\in\cC(\rho_m^*)\). Set \(\rho'=\rho_m^*\). This
satisfies~\eqref{eq:extension-conclusion}. We complete its domain in Step~5.

\smallskip
\noindent\textbf{Step 3. The answer is positive.}
Suppose instead that some oracle gives output \(1\). Choose
\(Y_m\in\cC(\rho_m^*)\) with \(\Gamma^{Y_m}(n)=1\), and an axiom
\(\tup{n,1,P_m,N_m}\in\Gamma\) witnessing the convergence.
Choose \(a_m\in M\) strictly above \(n\) and every element of
\(\dom\rho_m\cup P_m\cup N_m\). Since a pair is at least as large as
either coordinate, every column index queried by \(P_m\) or \(N_m\)
is below \(a_m\). In particular, for every \(k\geq a_m\),
\begin{equation}\label{eq:unused-columns}
 P_m^{[2k+1]}=P_m^{[2k+2]}=N_m^{[2k+1]}=N_m^{[2k+2]}=\varnothing.
\end{equation}
Let \(\rho_{m+1}\) have domain \([0,a_m)^M\setminus\{n\}\) and
record on this domain the column pairs of \(Y_m\). Since
\(Y_m\in\cC(\rho_m^*)\subseteq\cC(\rho_m)\), this extends
\(\rho_m\), and \(Y_m\in\cC(\rho_{m+1})\).

\smallskip
\noindent\textbf{Step 4. The positive answer occurs at every standard stage.}
Suppose Step~2 never applies. Let
\(\rho_\omega=\bigcup_{m\in\omega}\rho_m\). Then \(n\notin\dom\rho_\omega\).
By Lemma~\ref{lem:main1}~\ref{it:like3}, \(\cf(M)>\aleph_0\), so the
countably many bounds of the domains \(\dom\rho_m\) have a common
upper bound in \(M\). Thus \(\rho_\omega\) is a bounded configuration.
For each \(m\in\omega\), write
\[
 S_m=(P_m\cup N_m)^{[2n+2]}\cap I
\]
for the part of column \(2n+2\) of \(P_m\cup N_m\) that lies in \(I\).
Since \(Y_m\in\cC(\rho_m^*)\), we have \(Y_m^{[2n+2]}=D\), so
Lemma~\ref{lem:externally-finite} makes \(S_m\) externally finite.
By Lemma~\ref{lem:main1}~\ref{it:like2}, \(\abs I=\kappa>\aleph_0\), so
choose
\[
 l\in I\setminus\bigcup_{m\in\omega}S_m
\]
and set \(\rho'=\rho_\omega\cup\{\tup{n,(0,l)}\}\).

Fix \(Z\in\cC(\rho')\). Since \(\rho'(n)=(0,l)\), we have
\(A_Z(n)=0\). It suffices to show that \(\Gamma^Z(n)\ne0\).
Suppose instead that \(\tup{n,0,P,N}\in\Gamma\) witnesses \(\Gamma^Z(n)=0\).
Fix \(m\in\omega\). The axiom \(\tup{n,1,P_m,N_m}\) was chosen in
Step~3, so Lemma~\ref{lem:opposite-outputs} says that the two axioms are incompatible, which gives
\begin{equation}\label{eq:query-intersection}
 P\cap N_m\ne\varnothing\quad\text{or}\quad P_m\cap N\ne\varnothing.
\end{equation}
Since \(P\subseteq Z\) and \(N_m\cap Y_m=\varnothing\), every element of
\(P\cap N_m\) lies in \(Z\setminus Y_m\). Likewise every element of
\(P_m\cap N\) lies in \(Y_m\setminus Z\). A witness
to~\eqref{eq:query-intersection} therefore lies at a position where
\(Z\) and \(Y_m\) differ.

Both oracles have column \(0\) equal to \(D\). Since
\(Z\in\cC(\rho')\subseteq\cC(\rho_{m+1})\) and \(Y_m\in\cC(\rho_{m+1})\),
their column pairs agree for every \(k<a_m\) with \(k\ne n\).
For \(k\geq a_m\), the column pairs of \(P_m\) and \(N_m\) are empty
by~\eqref{eq:unused-columns}. The witness is therefore in column
\(2n+1\) or \(2n+2\), so at least one of the following holds.
\begin{align}
 P^{[2n+1]}\cap N_m^{[2n+1]}&\ne\varnothing,\label{eq:incomp-i}\\
 P_m^{[2n+1]}\cap N^{[2n+1]}&\ne\varnothing,\label{eq:incomp-ii}\\
 P^{[2n+2]}\cap N_m^{[2n+2]}&\ne\varnothing,\label{eq:incomp-iii}\\
 P_m^{[2n+2]}\cap N^{[2n+2]}&\ne\varnothing.\label{eq:incomp-iv}
\end{align}
The two pairs of columns are
\[
 (Z^{[2n+1]},Z^{[2n+2]})=(D,D\Delta\{l\}),\qquad
 (Y_m^{[2n+1]},Y_m^{[2n+2]})=(D\Delta\{m\},D).
\]
They differ in column \(2n+1\) only at \(m\), and in column \(2n+2\)
only at \(l\). Since \(l\in I\setminus S_m\), we have
\(l\notin P_m^{[2n+2]}\cup N_m^{[2n+2]}\), so neither
\eqref{eq:incomp-iii} nor~\eqref{eq:incomp-iv} holds. Thus the witness
is \(\tup{2n+1,m}\), and~\eqref{eq:incomp-i} or~\eqref{eq:incomp-ii}
gives \(m\in P^{[2n+1]}\cup N^{[2n+1]}\). Every standard \(m\)
lies in \(I\), so
\[
 \omega\subseteq P^{[2n+1]}\cup(N^{[2n+1]}\cap I).
\]
But \(Z^{[2n+1]}=D\), so Lemma~\ref{lem:externally-finite} makes both
sets on the right externally finite. This contradiction proves~\eqref{eq:extension-conclusion}.

\smallskip
\noindent\textbf{Step 5. Extending the domain to an initial segment.}
In either Step~2 or Step~4, \(\rho'\) is bounded and \(n\in\dom\rho'\).
Choose \(a\in M\) strictly above every element of \(\dom\rho'\).
Extend \(\rho'\) to a configuration \(\sigma\) on \([0,a)^M\) by assigning
\((0,0)\) to each unassigned \(k<a\). Since \(\cC(\sigma)\subseteq\cC(\rho')\),
\eqref{eq:extension-conclusion} holds for every \(Z\in\cC(\sigma)\).
\end{proof}

\begin{proof}[Proof of Lemma~\ref{lem:construction}]
Since \(\abs M=\kappa^+\) by Lemma~\ref{lem:main1}~\ref{it:like3},
list its elements as
\(\{e_\alpha\mid\alpha<\kappa^+\}\). We construct an increasing
sequence of bounded configurations \(\{\rho_\alpha\mid\alpha<\kappa^+\}\)
with initial segment domains.
For each \(\alpha\), we meet the requirement
\begin{equation}\label{eq:requirements}
 R_\alpha:\quad\exists n\,\bigl(\Phi_{e_\alpha}^X(n)\ne A_X(n)\bigr).
\end{equation}
The indices are elements of \(M\), including its nonstandard elements;
the order of construction is external.

Let \(\rho_0=\varnothing\). Given \(\rho_\alpha\),
choose \(n\) above \(\dom\rho_\alpha\) and apply
Lemma~\ref{lem:extension} to \(\Phi_{e_\alpha}\), \(\rho_\alpha\),
and \(n\). For the resulting \(\rho_{\alpha+1}\), every
\(Z\in\cC(\rho_{\alpha+1})\) satisfies
\[
 \Phi_{e_\alpha}^Z(n)\ne A_Z(n).
\]

At a nonzero limit \(\alpha<\kappa^+\), let
\[
 \rho_\alpha=\bigcup_{\beta<\alpha}\rho_\beta.
\]
Choose a bound for each \(\dom\rho_\beta\), \(\beta<\alpha\).
There are fewer than \(\kappa^+\) such bounds, and
\(\cf(M)=\kappa^+\) by Lemma~\ref{lem:main1}~\ref{it:like3}, so they
have a common upper bound in \(M\). Thus \(\rho_\alpha\) is again a
bounded configuration.

Finally, let \(\rho_\infty=\bigcup_{\alpha<\kappa^+}\rho_\alpha\).
Each successor stage adds a new input to the domain, so
\(\abs{\dom\rho_\infty}=\kappa^+\). Since \(\dom\rho_\infty\) is an
initial segment, Lemma~\ref{lem:main1}~\ref{it:like2} gives
\(\dom\rho_\infty=M\). Let \(X\) be the unique element of
\(\cC(\rho_\infty)\). Then \(X\in\cC(\rho_{\alpha+1})\) for every
\(\alpha<\kappa^+\). Every Turing functional occurs in our list, so
all requirements~\eqref{eq:requirements} are met.
Therefore \(A_X\notin\pT(X)\), while Lemma~\ref{lem:enumerations}
gives \(A_X\in\CEpm(X)\).

We emphasize that the configurations and the first extension lemma are
external objects and assertions. This recursion is carried out in ZFC;
it is not an application of induction inside \(M\).
\end{proof}

\subsection{Passing to a countable model}\label{sec:countable}

\begin{lemma}\label{lem:countable_model}
Let \(M\models I\Sigma_1\) and \(A,X\subseteq M\) satisfy
\(A\in\CEpm(X)\setminus\pT(X)\). There exists a countable model
\(M'\) and \(A',X'\subseteq M'\) such that
\((M',A',X')\prec(M,A,X)\) and
\(A'\in\CEpm(X')\setminus\pT(X')\) in \(M'\).
\end{lemma}
\begin{proof}
Expand the language by unary predicates for \(A\) and \(X\).
Using the fixed enumeration of Turing functionals, each of the assertions
that \(A\) is \(X\)-c.e., that \(\overline{A}\) is \(X\)-c.e., and that
\(A\notin\pT(X)\) is first-order expressible in this language.
Apply downward L\"owenheim--Skolem and use elementarity.
\end{proof}

Take \(\kappa=2^{\aleph_0}\) and the model \(M\) of
Lemma~\ref{lem:main1}. By Lemma~\ref{lem:construction}, there are
\(A,X\subseteq M\) with \(A\in\CEpm(X)\setminus\pT(X)\).
Lemma~\ref{lem:countable_model} gives a countable \(N\equiv M\equiv
\mathbb N\) in which this separation still holds. This proves
Theorem~\ref{thm:main}(\ref{it:main-pointwise}).

\subsection{All three separations for the same oracle}\label{sec:other-inclusions}

We now return to the large model \(M\), keeping \(b\), \(I=[0,b)^M\),
and \(D\) from Section~\ref{sec:construction}. The oracle constructed
there already separates strong from pointwise Turing reducibility.
To see this, for each column oracle \(Y\) let
\begin{equation}\label{eq:column-agreement}
 B_Y=\{\tup{n,k}\mid Y^{[2n+1]}(k)=Y^{[0]}(k)\}.
\end{equation}
\begin{lemma}\label{lem:agreement-pointwise}
For every column oracle \(Y\), \(B_Y\le_{pT}Y\).
\end{lemma}
\begin{proof}
Given \(\tup{n,k}\), let \(a_{n,k}=\tup{2n+1,k}\) and \(b_k=\tup{0,k}\).
To decide whether \(\tup{n,k}\in B_Y\), we compare the two bits of the
oracle at \(a_{n,k}\) and \(b_k\). More precisely, for \(P\cap N=\varnothing\),
we enumerate \(\tup{\tup{n,k},1,P,N}\) into a Turing functional \(\Phi\)
if either \(a_{n,k},b_k\in P\) or \(a_{n,k},b_k\in N\). We enumerate
\(\tup{\tup{n,k},0,P,N}\) if one of \(a_{n,k},b_k\) belongs to \(P\)
and the other to \(N\). Then \(B_Y=\Phi^Y\).
\end{proof}

As in Definition~\ref{def:turing_reducibility}, \(B_Y^+=(B_Y)^+\)
is the set of codes of \(M\)-finite subsets of \(B_Y\). Membership in
\(B_Y^+\) is defined by a bounded formula with predicate \(Y\). For an \(M\)-finite set
\(F\) with code \(c\),
\begin{equation}\label{eq:B-Y-plus}
 c\in B_Y^+\iff
 \forall n,k\le c\,(\tup{n,k}\in F\to Y^{[2n+1]}(k)=Y^{[0]}(k)).
\end{equation}

\begin{lemma}\label{lem:pT-not-sT}
Let \(X\) be a column oracle with \(A_X\notin\pT(X)\). Then \(B_X\in\pT(X)\setminus\sT(X)\).
\end{lemma}
\begin{proof}
    By Lemma~\ref{lem:agreement-pointwise}, \(B_X\in\pT(X)\).
    For \(n\in M\), let \(f(n)\) be the code of the \(M\)-finite set
    \(\{n\}\times I=\{\tup{n,k}\mid k<b\}\). Then \(f\) is computable
    with parameter \(b\). By~\eqref{eq:A-Y},
    \[
     n\notin A_X\iff f(n)\in B_X^+.
    \]
    If \(B_X\le_{sT}X\), then \(B_X^+=\Theta^X\) for some Turing functional
    \(\Theta\). On input \(n\), we can compute \(f(n)\) and then run
    \(\Theta^X(f(n))\) to obtain \(\overline{A_X}(n)\). Reversing the
    output gives \(A_X\le_{pT}X\), a contradiction.
\end{proof}

Thus the oracle witnessing part~(\ref{it:main-pointwise}) in the large model
also witnesses part~(\ref{it:main-strong}).

For Theorem~\ref{thm:main}(\ref{it:main-definability}), \eqref{eq:B-Y-plus} already gives
\(B_X^+\in\Delta_0^X\subseteq \Delta_1^X\). We add one requirement for each Turing
functional \(\Gamma\) to ensure that \(B_X^+\ne\dom\Gamma^X\).
The following lemma lets us meet this requirement by extending the
configuration obtained from Lemma~\ref{lem:extension}.

\begin{lemma}[Second extension lemma]\label{lem:second-extension}
Let \(\rho\) be a bounded configuration, let \(\Gamma\) be a Turing
functional, and let \(t\in M\) be above \(\dom\rho\). There exist
\(g\in M\) and a bounded configuration \(\rho'\supseteq\rho\) such
that \(\dom\rho'\) is an initial segment containing \(t\) and
\[
 g\in B_Z^+\mathbin{\Delta}\dom\Gamma^Z
 \qquad\text{for every }Z\in\cC(\rho').
\]
\end{lemma}
\begin{proof}
Let \(g\) be the code of \(G=\{\tup{t+j,j}\mid j<b\}\).
Extend \(\rho\) to a configuration \(\sigma\) with domain \([0,t+b)^M\)
by putting \(\sigma(t+j)=(0,j)\) for \(j<b\) and assigning \((0,0)\)
to each remaining unassigned input below \(t\). For \(Y\in\cC(\sigma)\)
and \(j<b\), we have
\begin{equation}\label{eq:second-extension-columns}
 (Y^{[2(t+j)+1]},Y^{[2(t+j)+2]})=(D,D\Delta\{j\}).
\end{equation}
Thus \(G\subseteq B_Y\), or equivalently, \(g\in B_Y^+\), for every such \(Y\). If there is no
\(Y\in\cC(\sigma)\) with \(\Gamma^Y(g)\downarrow\), take
\(\rho'=\sigma\).

Otherwise choose such a \(Y\) and an axiom \(\tup{g,y,P,N}\in\Gamma\)
witnessing convergence. For \(j<b\), put
\[
 u_j=\tup{2(t+j)+1,j},\qquad v_j=\tup{2(t+j)+2,j}.
\]
By~\eqref{eq:second-extension-columns}, if \(u_j\in P\subseteq Y\), then \(D(j)=Y(u_j)=1\), and if \(v_j\in N\subseteq \overline{Y}\), then \(D(j)=1-Y(v_j)=1\). Hence
\[
 \{j<b\mid u_j\in P \ \lor\ v_j\in N\}\subseteq D.
\]
Similarly,
\[
    \{j<b\mid u_j\in N \ \lor\ v_j\in P\}\subseteq I\setminus D.
\]
Both sets are \(M\)-finite, so
both are externally finite by Lemma~\ref{lem:D}. Choose \(j\in I\)
outside their union. In the pair of columns for \(t+j\), replace
\((D,D\Delta\{j\})\) by \((D\Delta\{j\},D)\), and call the
resulting oracle \(Y'\). Only the bits at \(u_j\) and \(v_j\) have
changed. Neither was queried, so the same axiom witnesses
\(\Gamma^{Y'}(g)\downarrow\). But \(\tup{t+j,j}\notin B_{Y'}\), and
hence \(g\notin B_{Y'}^+\).

Choose \(c\ge t+b\) above every element of \(P\cup N\) and let
\(\rho'\) assign the column pairs below \(c\) according to \(Y'\) so that \(Y'\in \cC(\rho')\).
Since \(t\) is above \(\dom\rho\), this extends \(\rho\). Every
\(Z\in\cC(\rho')\) agrees with \(Y'\) on \(P\cup N\) and on the
column pair for \(t+j\). Thus \(\Gamma^Z(g)\downarrow\) and
\(g\notin B_Z^+\), as required.
\end{proof}

\begin{proof}[Proof of Theorem~\ref{thm:main}]
In the model \(M\) of Lemma~\ref{lem:main1}, carry out the construction
of Lemma~\ref{lem:construction} with one more extension at each successor
stage. For \(\Phi_{e_\alpha}\), first apply Lemma~\ref{lem:extension}
with \(n\) above the current domain. Then apply
Lemma~\ref{lem:second-extension} to the resulting configuration and the
same functional, choosing \(t\) above the new domain. The first
requirement continues to hold after the second extension, since
\(\cC(\rho')\subseteq\cC(\rho)\) whenever \(\rho'\supseteq\rho\). At limit stages take unions as before.

For the resulting column oracle \(X\), the first extension lemma gives
\(A_X\notin\pT(X)\), and the second gives
\(B_X^+\ne\dom\Phi_e^X\) for every \(e\in M\). By
Lemma~\ref{lem:enumerations}, \(A_X\in\CEpm(X)\). Thus
Lemma~\ref{lem:pT-not-sT} and~\eqref{eq:B-Y-plus} give
\[
 \sT(X)\subsetneq\pT(X)\subsetneq\CEpm(X)\subsetneq\Delta_1^X,
\]
with witnesses \(B_X\), \(A_X\), and \(B_X^+\), respectively.

As in Section~\ref{sec:countable}, take a countable elementary
substructure of \((M,A_X,X)\) containing \(b\). The definitions
\eqref{eq:column-agreement} and~\eqref{eq:B-Y-plus} define \(B_X\)
and \(B_X^+\) from \(X\), so the assertions about all three witnesses
pass to this substructure by elementarity. Its arithmetic reduct
\(N\) satisfies \(N\equiv M\equiv\mathbb N\), proving the theorem.
\end{proof}

We give a final remark that the set \(D\) is not definable in our model \(M\) as in a model of PA, every definable subset of \(M\) is regular, this contradicts the fact that \(D\) cannot be regular (Lemma~\ref{lem:D}). Therefore, \(X\) is also not definable in \(M\).

\section*{Acknowledgements}
The question studied in this paper was formulated by the authors. The idea for the proof of Theorem~\ref{thm:main} was suggested by GPT-6 Astra. Based on this idea, the authors prepared the initial manuscript. Subsequently, AI tools were used to check the arguments and improve the exposition.
The authors remain responsible for the content of the paper and for any errors. We also thank Yue Yang for helpful comments on the manuscript.

\bibliographystyle{plainnat}
\bibliography{ref}
\end{document}